\documentclass[12pt]{amsart}
\usepackage{epsfig}
\usepackage{amsmath,amssymb}
\usepackage{youngtab}
\numberwithin{equation}{section}

\newtheorem{theorem}{Theorem}[section]

\newtheorem{remark}[theorem]{Remark}

\allowdisplaybreaks \numberwithin{equation}{section}

\begin{document}
\title[$\sigma$- and $\tau$-functions]{Sigma, tau and Abelian
  functions of algebraic curves }

\author{J. C. Eilbeck} \address{Department of Mathematics and the
  Maxwell Institute for Mathematical Sciences, Heriot-Watt University,
  Edinburgh, UK EH14 4AS} \email{J.C.Eilbeck@hw.ac.uk}
\author{V.Z. Enolski} \address{ Hanse-Wissenschaftskolleg, 27753 Delmenhorst, Germany,
\\On leave from Institute of Magnetism, National
  Academy of Sciences of Ukraine, Kyiv-142, Ukraine }
\email{vze@ma.hw.ac.uk}

\author{J.Gibbons}
\address{Imperial College, 180 Queen's Gate, London SW7 2BZ}
\email{j.gibbons@ic.ac.uk}
\date{\today}

\begin{abstract}
  We compare and contrast three different methods for the construction
  of the differential relations satisfied by the fundamental Abelian
  functions associated with an algebraic curve.  We realize these
  Abelian functions as logarithmic derivatives of the associated sigma
  function.  In two of the methods, the use of the tau function,
  expressed in terms of the sigma function, is central to the
  construction of differential relations between the Abelian
  functions.

\end{abstract}

\maketitle


\section{Introduction}

In the mid-1970s, the results of Its, Matveev, Dubrovin and Novikov (see
\cite{dmn76}) led to the discovery of a remarkable $\theta$-functional
formula to solve the KdV equation $u_t=6uu_x-u_{xxx}$. This solution
was given as a second logarithmic derivative of a Riemann
theta-function:
\begin{equation}
u(x,t)=-\frac{\partial^2}{\partial x^2 }  \,\mathrm{ln} \,
\theta(\boldsymbol{U}x+\boldsymbol{V}t+ \boldsymbol{W})+C
\end{equation}
with $\boldsymbol{U}, \boldsymbol{V}, \boldsymbol{W} = \mathrm{const}
\in \mathbb{C}^g$ and $C\in \mathbb{C}$. The theta-function in that
formula was constructed from a hyperelliptic curve $X_g$ of genus $g$,
while the ``winding vectors'' $\boldsymbol{U}, \boldsymbol{V}$ are
periods of abelian differentials of the second kind on $X_g$. Further,
this formula is in a sense universal; it was generalized by Krichever
\cite{kr77} to other integrable hierarchies - whose solutions were
associated with other algebraic curves.  In this paper we consider the
converse problem:

{\em Given an algebraic curve $X_g$, of genus $g$, its Riemann period
  matrix $\tau$, and its Jacobi variety $\mathrm{Jac}(X_g) =
  \mathbb{C}/(1_g\oplus \tau)$, we may construct $\theta$-functions
  $\theta(\boldsymbol{u};\tau)$, $\boldsymbol{u}\in \mathrm{Jac}(X_g)
  $; then the fundamental Abelian functions on $\mathrm{Jac}(X_g)$ may
  be realized as the second logarithmic derivatives of
  $\theta(\boldsymbol{u},\tau)$, $\wp_{ij}(\boldsymbol{u}) = -
  \frac{\partial^2\ln \theta(\boldsymbol{u};\tau) }{\partial u_i
    \partial u_j}$.  We wish to construct all differential relations
  between these Abelian functions on $\mathrm{Jac}(X_g) $.}

The simplest case, the Weierstrass cubic, $y^2=4x^3-g_2x-g_3$, is an
algebraic curve of genus one.  This is uniformized by the Weierstrass
elliptic functions, $x=\wp(u), y=\wp'(u)$, and the differential
relations are
\begin{equation}
\wp''=6\wp^2-\frac{g_2}{2},\quad {\wp'}^2=4\wp^3-g_2\wp-g_3.
\label{weierstrass}
\end{equation}
In the case of higher genera, $g>1$, the derivation of analogous
equations becomes much more complicated.  In particular, the
fundamental Abelian functions are now {\em partial} derivatives of a
function of $g$ variables.  The different approaches to this problem
form the main content of the paper.  We restrict our analysis to the
case of $(n,s)$-curves introduced and investigated in this context by
Buchstaber, Leykin and Enolski \cite{bel97b,bel99}
\begin{equation}
X_g:\;\;y^n=x^s+\sum_{ni+sj<ns} \lambda_{ij} x^iy^j.
\end{equation}
These represent a natural generalization of elliptic curves to higher
genera, and include the general hyperelliptic curve ($n=2$).

To any such curve we may associate an object which is fundamental to
all our treatments of this problem, the fundamental bi-differential.
This is the unique symmetric meromorphic $2$-form on $X_g\times X_g$,
whose only second order pole lies on the diagonal $Q=S$, and which
satisfies
\begin{align*}
\omega(Q,S)-\frac{\mathrm{d}\xi(Q)\mathrm{d}\xi(S)}{(\xi(Q)-\xi(S))^2 }
=\Phi(\xi(S),\xi(Q))\mathrm{d}\xi(Q)\mathrm{d}\xi(S),
\end{align*}
where $\Phi(\xi(S),\xi(Q))$ is holomorphic, and $\xi(Q),\xi(S)$ are
local coordinates in the vicinity of a base point $P$, $\xi(P)=0$.
Usually $\omega(Q,S)$ is realized as the second logarithmic derivative
of the prime-form or theta-function \cite{fa73}.  But in our
development we use an alternative representation of $\omega(Q,S)$ in
the {\em algebraic form} that goes back to Weierstrass and Klein, and
which was well documented by Baker \cite{ba97}
\begin{equation}
  \omega(Q,S)=\frac{\mathcal{F}(Q,S)}{f_y(Q) f_w(S)(x-z)^2}
  \mathrm{d}x\mathrm{d}z+
  2\mathrm{d}\boldsymbol{u}(Q)^T\varkappa \mathrm{d}
  \boldsymbol{u}(S), \label{omegaQS}
\end{equation}
where $Q=(x,y)$, $S=(z,w)$, and the function
$\mathcal{F}(Q,S)=\mathcal{F}((x,y),(z,w))$ is a polynomial of its
arguments with coefficients depending on the parameters of the curve
$X_g$ $\mathrm{d}\boldsymbol{u}=(\mathrm{d}u_1,\ldots,\mathrm{d}u_g )^T$ is the vector of basic holomorphic differentials. 
The factor 2 in the final term is introduced so that this
theory reduces to the Weierstrass theory in the elliptic case.  Here,
$\varkappa$ is a symmetric matrix expressed in terms of the first and
second period matrices, $2\omega$, $2\eta$ respectively, as
$\varkappa= \omega^{-1}\eta$.  This provides the normalization of
$\omega(Q,S)$.  We will refer to the first term on the right of
(\ref{omegaQS}), which involves the polynomial $\mathcal{F}(Q,S)$, as
the algebraic part, so that
\begin{equation}
  \omega(Q,S)=\omega^{\rm{alg}}(Q,S)+ 2\mathrm{d}
  \boldsymbol{u}(Q)^T\varkappa \mathrm{d} \boldsymbol{u}(S).
\end{equation}
This representation was revisited and developed by Buchstaber, Leykin
and Enolski \cite{bel97b} and more recently by Nakayashiki
\cite{nakaya08}.

The algebraic representation of the fundamental differential, as
described above, lies behind the definition of the multivariate
sigma function in terms of the theta-function. This differs from $\theta$
by an exponential factor and a modular factor:
\begin{equation}
  \sigma(\boldsymbol{u}) = C(\tau)\mathrm{exp} \left\{
    \tfrac12 \boldsymbol{u}^T \omega^{-1}\eta \boldsymbol{u}  \right\}
  \theta\left(\tfrac12 \omega^{-1}\boldsymbol{u};\tau \right). \label{sigma}
\end{equation}
Here the $g\times g$ matrices $2\omega, 2\eta$ are the first and
second period matrices, and $\tau=\omega^{-1}\omega'$.  The modular
constant $C(\tau)$ is known explicitly for hyperelliptic curves and a
number of other cases. However its explicit form is not necessary here, for
the fundamental Abelian functions are independent of $C(\tau)$.
These modifications make $\sigma(\boldsymbol{u})$ invariant with
respect to the action of
the symplectic group, so that for any $\gamma\in{\mathrm Sp}(2g,
\mathbb{Z})$, we have:
\begin{equation}
  \sigma(\boldsymbol{u};\gamma \tau)=\sigma(\boldsymbol{u};\tau).
\end{equation}
The multivariate sigma-function is the natural generalization of the
Weierstrass sigma function to algebraic curves of higher genera,
i.e.\ $(n,s)$-curves in this context. In his lectures \cite{w-lect},
Weierstrass started by defining the sigma-function in terms of series
with coefficients given recursively, which was the key point of the
Weierstrass theory of elliptic functions. A generalization of this
result to the genus two curve was started by Baker \cite{ba07} and
recently completed by Buchstaber and Leykin \cite{bl05}, who obtained
recurrence relations between coefficients of the sigma-series in
closed form.  In addition, Buchstaber and Leykin recently found an
operator algebra that annihilates the sigma-function of a higher
genera $(n,s)$-curve \cite{bl08}. The recursive definition of the
higher genera sigma-functions remains a challenging problem to solve,
with \cite{bl08} providing a definite step. We believe that the future
theory of the sigma and corresponding Abelian functions can be formulated
on the basis of sigma expansions that will complete the extension of
the Weierstrass theory to curves of higher genera.

In this paper we study the interrelation of the multivariate sigma and
Sato tau functions. The $\tau$-function was introduced by Sato
\cite{sato80,sato81} in the much more general context of integrable
hierarchies. But it seems there are few results studying algebraic
curves in Sato theory. However we should mention recent work by
Konopelchenko and Ortenzi analyzing algebro-geometric structure in
Birkhoff strata of the Sato Grassmannian \cite{ko10}.  The recent
papers of Matsutani and Previato \cite{matprev08,matprev10}, studying
Jacobi inversion on Jacobian strata of $(r,s)$ curves, is also
relevant to our work, relating stratification of the Sato Grassmannian
to partitions.

Here we deal with the `algebro-geometric
$\tau$-function' (AGT) associated with an algebraic curve.

The AGT of the genus $g$ curve $X_g$ is defined, following Fay,
\cite{fay83,fay89}, as a function of the `times'
$\boldsymbol{t}=(t_1,\ldots,t_g,t_{g+1},\ldots)$, a point
$\boldsymbol{u}\in\mathrm{Jac}(X_g)$, as well as a point $P\in X_g$;
it is given by the formula
\begin{align*}
\tau(\boldsymbol{t};\boldsymbol{u},P) = \theta\left( \sum_{k=1}^{\infty}
\boldsymbol{U}_k(P) t_k+\frac12\omega^{-1}\boldsymbol{u} \right) \mathrm{exp}
\left\{ \frac12\sum_{m,n\geq 1} \omega_{mn}(P) t_m t_n \right\}.
\end{align*}
Here the ``winding vectors'' $\boldsymbol{U}_k(P)$ appear in the
expansion of the normalized holomorphic integral $\boldsymbol{v}$, the
quantities $\omega_{mn}(P)$ define the holomorphic part of the
expansion of the fundamental differential of the second kind
$\omega(Q,S)$ near the point $P$.
We then introduce the $\tau$-function by the
formula:
\begin{align}
  \frac{ \tau(\boldsymbol{t};\boldsymbol{u},P)}
  {\tau(\boldsymbol{0};\boldsymbol{u},P)}= \frac{
    \sigma\left(\sum_{k=1}^{\infty} \mathcal{A}^{-1}\boldsymbol{U}_k(P) t_k+
      \boldsymbol{u}\right)} {\sigma(\boldsymbol{u})}\mathrm{exp}
  \left\{ \frac12 \sum_{k,l=0}^{\infty}
    \omega^{\mathrm{alg}}_{k,l}(P)t_k t_l \right\}.
\label{tausigma}\end{align}

This representation of $\tau$ in terms of $\sigma$ was used by Harnad
and Enolski \cite{enhar08} to analyze the Schur function expansion of
$\tau$ for the case of algebraic curves. Recently Nakayashiki
\cite{nakayashiki09} has independently suggested a similar expression
for the AGT in terms of multivariate $\sigma$-functions. In this paper
we concentrate on the application of this representation to the
derivation of the differential relations between Abelian functions of
the $(n,s)$-curve, continuing and developing the work of
\cite{enhar08}.

Developing a further analogy with the Weierstrass theory of elliptic
functions, we represent the Abelian
functions, that is, $2g$-periodic
functions
\begin{align*}
F(\boldsymbol{u}+2\boldsymbol{n}\omega +2\boldsymbol{n}'\omega')
  & = F(\boldsymbol{u}), \qquad \forall \boldsymbol{n}, \boldsymbol{n}' \in
  \mathbb{N}
  \end{align*}
  on the Jacobian
  \begin{align*}
  \widetilde{\mathrm{Jac}}(X_g) & =\mathbb{C}/2\omega \oplus 2\omega'=
  \mathcal{A}^{-1}\mathrm{Jac}(X_g),\quad 2\omega =\mathcal{A},
\end{align*}
as second and higher logarithmic derivatives:
\begin{align}
\zeta_i(\boldsymbol{u})&=\frac{\partial}{\partial u_i}
\mathrm{ln}\,\sigma(\boldsymbol{u}),\label{zeta} \\
\wp_{ij}(\boldsymbol{u})&=-\frac{\partial^2}{\partial u_i \partial u_j}
\mathrm{ln}\,\sigma(\boldsymbol{u}), \quad \wp_{ijk}(\boldsymbol{u})=
-\frac{\partial^3}{\partial u_i \partial u_j \partial u_k }
\mathrm{ln}\,\sigma(\boldsymbol{u}),\quad \rm{etc}., \label{Kleinwp}
\end{align}
where $i,j,k, \dots = 1\ldots,g$. We should remark that the
$\zeta_i(\boldsymbol{u})$ are not Abelian functions.  In this
notation, the genus one Weierstrass equations (\ref{weierstrass}) become
 \[
 \wp_{1111}=6\wp_{11}^2-\frac{g_2}{2},\quad {\wp_{111}}^2 =
 4\wp_{11}^3-g_2\wp_{11}-g_3.
\]
For general $g$, the $\wp_{ij}, \wp_{ijk}, \ldots$ are called Kleinian $\wp$-functions. They
are convenient coordinates to represent the dependent variables in
the hierarchy of integrable systems.

In this paper we compare and contrast three approaches to obtain the
partial differential relations for the Abelian functions associated
with the $(n,s)$-curve $X_g$, using Kleinian $\wp$-functions as
coordinates.

The first of these, and the best known, is the classical approach of
comparing two different expansions of the fundamental bi-differential;
this yields first the solution of the Jacobi inversion problem for the
curve, and in higher orders, a sequence of differential relations
involving the $\wp_{ij}$.

The $\tau$-function approach to the derivation of completely
integrable systems of KP type has led to two different ways
\cite{djkm83} to obtain relations between Taylor coefficients of the
$\tau$ function expansion.  The first of these specifically exploits
the fact that these Taylor coefficients are determinants, i.e.\
Pl\"ucker coordinates in the Grassmanian, and they hence satisfy the
Pl\"ucker relations \cite{fay83}. The second is based on the Bilinear
Identity, which leads to the Residue Formula \cite{fay89}, giving a
family of differential polynomials in $\tau$ which must vanish - these
are the partial differential equations we need.

We further consider and compare the two techniques based on the
$\tau$-function method, which give a derivation of the required
differential relations; specializing to a particular algebraic curve,
we consider its algebro-geometric $\tau$-function.  The special
feature of our development is that we define this $\tau$-function in
terms of the multidimensional $\sigma$-function of the curve, leading
to coordinates that are explicitly written in terms of Kleinian
$\wp$-functions.  The differential relations we find between these
functions can be understood as arising from special solutions of
integrable hierarchies of KP type, associated with the given curve
(see for example \cite{ba97,bel97b,nakaya08}).  We will describe the
correspondence between individual differential equations for
$\wp$-functions with Young diagrams defining Pl\"ucker relations.  We
illustrate these approaches by considering two particular examples:
the genus 2 hyperelliptic curves \cite{ba07}, and the genus 3 trigonal
curve, which can be found in \cite{bel97b,bel00,eemop07}).

We remark that in similar context, Nimmo \cite{nimmo} in 1982 applied
the methods of symmetric function theory to describe the KP hierarchy,
and studied the action of the recursion operator.  The general
approach based on Pl\"ucker coordinates for deriving KP-flows in terms
of Kleinian $\sigma$-functions was recently discussed in
\cite{enhar08}. Here we develop these ideas and consider some
non-trivial examples to clarify the interrelation of the
$\tau$-functional formulation of integrable hierarchies and the
$\sigma$-functional approach.  We will also consider the relationship
between of this derivation, based on the Pl\"ucker relations, and that
based on the Residue Formula.  Both give a systematic way of
generating the required relations, but the differences between the two
approaches are instructive.

\section{Algebraic curves}
Let $X_g$ be a genus $g\geq 1$ algebraic curve given by the polynomial
equation
\begin{equation}
f(x,y)=0,\quad f(x,y)=y^n+y^{n-1}a_1(x)+\ldots+a_0(x).\label{gencurve}
\end{equation}
We shall consider in what follows two relatively simple curves of the
class (\ref{gencurve}),

{\sf  Example I}:  the hyperelliptic genus two curve
\begin{equation}
y^2=4x^5+\alpha_4 x^4+\ldots+\alpha_0\label{heperel}
\end{equation}
and

{\sf  Example II}:   the cyclic trigonal genus three curve

\begin{equation}
f(x,y)=y^3 -(x^4 +\mu_3 x^3 +\mu_6 x^2 +\mu_9 x +\mu_{12}).\label{trig}
\end{equation}

We equip $X_g$ with a canonical basis of cycles
$(\mathfrak{a}_1,\ldots,\mathfrak{a}_g;\mathfrak{b}_1,
\ldots,\mathfrak{b}_g)\in H_1(X,\mathbb{Z})$.  We denote by
$\mathrm{d}\boldsymbol{u}=(\mathrm{d} u_1,\ldots,\mathrm{d}u_g)^T$ the
vector whose entries are independent holomorphic differentials of the
curve $X_g$ as well as their $\mathfrak{a}$ and $\mathfrak{b}$-periods,
\begin{equation}
  2\omega=\left(\oint_{\mathfrak{a}_j} \mathrm{d}u_i\right)_{i,j=1,\ldots,g},
  \quad  2\omega'=\left(\oint_{\mathfrak{b}_j}
    \mathrm{d}u_i\right)_{i,j=1,\ldots,g}
\end{equation}
The period matrix $(2\omega,2\omega')$ is the first period matrix, and the
matrix $\tau=\omega^{-1}\omega'$ belongs to the upper Siegel
half-space, $\mathfrak{S}: \tau^T=\tau,\; \mathrm{Im}\, \tau >0 $.

The $\theta$-function $\theta[\alpha](\boldsymbol{z};\tau)$ with
characteristics $ [\alpha]=\left[ \begin{array}{c}
    \boldsymbol{a}^T\\\boldsymbol{b}^T
  \end{array} \right] $ , $[2\alpha]\in \mathbb{Z}^g\times\mathbb{Z}^g $
of the algebraic curve $X_g$ of genus $g$ is defined through its Fourier series
\begin{equation*}
  \theta\left[\begin{array}{c}\boldsymbol{a}^T\\\boldsymbol{b}^T
\end{array}\right](\boldsymbol{z};\tau)
  =\sum_{\boldsymbol{m}\in\mathbb{Z}^g}
  \mathrm{exp}\left\{ \imath\pi \left( \boldsymbol{m}+
      \boldsymbol{a}\right)^T\tau \left( \boldsymbol{m}+
      \boldsymbol{a}\right)+ 2\imath\pi\left( \boldsymbol{m}
      +\boldsymbol{a}\right)^T\left( \boldsymbol{z}
      +\boldsymbol{b}\right)
  \right\}
\end{equation*}
and possesses the periodicity property: for arbitrary
$\boldsymbol{a},\boldsymbol{b}\in \mathbb{Q}^g$ and arbitrary
$\boldsymbol{a}',\boldsymbol{b}'\in \mathbb{Q}^g$ the following
formula is valid
\begin{align*}
\begin{split}
&\theta\left[\begin{array}{c}\boldsymbol{a}^T\\\boldsymbol{b}^T
  \end{array}\right] (\boldsymbol{z}+\tau
\boldsymbol{a}'+\boldsymbol{b}' ;\tau)\\
&=\theta\left[\begin{array}{c}(\boldsymbol{a}+\boldsymbol{a}')^T
    (\boldsymbol{b}+\boldsymbol{b}')^T
  \end{array}\right](\boldsymbol{z} ;\tau) \mathrm{exp}\left\{
  -\imath\pi {\boldsymbol{a}'}^T\tau {\boldsymbol{a}'} -2\imath\pi
  {\boldsymbol{a}'}^T\boldsymbol{z} -2\imath\pi
  (\boldsymbol{b}+{\boldsymbol{b}'})^T\boldsymbol{a}' \right\}
\end{split}
\end{align*}

We introduce the associated meromorphic differentials
$\mathrm{d}\boldsymbol{r}=(\mathrm{d} r_1,\ldots,\mathrm{d}r_g)^T$ and
their periods
\begin{equation*}
  2\eta=-\left(\oint_{\mathfrak{a}_j} \mathrm{d}r_i\right)_{i,j=1,\ldots,g},
\quad  2\eta'=-\left(\oint_{\mathfrak{b}_j} \mathrm{d}r_i\right)_{i,j=1,\ldots,g}
\end{equation*}
which form the second period matrix  $(2\eta,2\eta')$.
The period matrices satisfy the condition
\begin{equation*}
  \left( \begin{array}{cc} \omega&\omega'\\ \eta&\eta' \end{array} \right)
  \left( \begin{array}{cc}0&1_g\\-1_g&0   \end{array}\right)
  \left( \begin{array}{cc} \omega&\omega'\\ \eta&\eta' \end{array} \right)^T
  =-\frac{\imath\pi}{2} \left( \begin{array}{cc}0&1_g\\-1_g&0
    \end{array}\right).
\end{equation*}
Here we denote the half-periods of the holomorphic and meromorphic
differentials by $(\omega,\omega')$ and $(\eta,\eta')$ in order to
emphasize the analogy with the Weierstrass theory. We will
also use the notation $\mathcal{A}=2\omega$ and $\mathcal{B}=2\omega'$
for the periods of holomorphic differentials. Further we denote
\begin{equation}
\mathrm{d}\boldsymbol{v}(Q)=(\mathrm{d} v_1(Q),\ldots,\mathrm{d}
v_g(Q))^T=\mathcal{A}^{-1} \mathrm{d}\boldsymbol{u}(Q) \label{vnormalized}
\end{equation}
as the vector of normalized holomorphic differentials.

The explicit calculation of canonical holomorphic differentials and
the meromorphic differentials conjugate to them is well understood;
in particular we have
\begin{align*}
  \mbox{\qquad{\sf Example I}}&
  \\
  \mathrm{d}u_1 & =\frac{x\,\mathrm{d}x}{y},\qquad \mathrm{d}u_2
=\frac{\mathrm{d}x}{y},\qquad \qquad \qquad \qquad \qquad \\
  \\
  \mathrm{d}r_1&=\frac{x^2\,\mathrm{d}x}{y},\qquad\mathrm{d}r_2
=\frac{x(\alpha_3+2 \alpha_4 x+12 x^2)\,\mathrm{d}x}{4y};\\
  \\
  \mbox{\qquad{\sf Example II}}&\\
\mathrm{d}u_1&=\frac{\mathrm{d}x}{3y},\qquad
\mathrm{d}u_2=\frac{x\mathrm{d}x}{3y^2},\quad\quad
\mathrm{d}u_3=\frac{\mathrm{d}x}{3y^2},\\
\\
\mathrm{d}r_1&=\frac{x^2\,\mathrm{d}x}{3y^2},\quad
\mathrm{d}r_2=-\frac{2x\,\mathrm{d}x}{3y},\quad
\mathrm{d}r_3=-\frac{(5x^2+3\mu_3x+\mu_6)}{3y}\,\mathrm{d}x.
\end{align*}

\begin{remark} Note that our labelling of the differentials is the
  reverse of \cite{eemop07}, with the interchange $1 \leftrightarrow
  2$ in {\sf example I}, and $(1,2,3) \leftrightarrow (3,2,1)$ in {\sf
    example II}.
\end{remark}

We introduce the fundamental bi-differential $\omega(Q,S)$ on $X\times X$
which is uniquely defined by the conditions:

{\bf(i)}  it is symmetric:
\begin{equation}
  \omega(Q,S)=\omega(S,Q);\label{omegasym}
\end{equation}

{\bf(ii)} it has its only pole along
the diagonal $Q=S$, in which neighbourhood it is expanded in a power
series according to
\begin{align}
  \omega(Q,S)-\frac{\mathrm{d}\xi(Q)\mathrm{d}\xi(S)}{(\xi(Q)-\xi(S))^2
  } =\sum_{m,n\geq 1} \omega_{mn}(P)
  \xi(Q)^{m-1}\xi(S)^{n-1}\mathrm{d}\xi(Q) \mathrm{d}\xi(S);
\label{omegapole}
\end{align}

{\bf(iii)}  it is normalized such that:
\begin{equation}
\oint_{\mathfrak{a}_j}\omega(Q,S)=0,\qquad j =1, \ldots, g.\label{omeganorm}
\end{equation}
The well known realization of the differential $\omega(Q,S)$ involves
the Schottky-Klein prime form $E(Q,S)$, which is a $(-1/2,-1/2)$-differential
defined for arbitrary points $Q,S\in X$
\begin{equation}
E(Q,S)=\frac{ \theta[\alpha]\left( \int_{Q}^S \mathrm{d} \boldsymbol{v}
\right) }{ h_{\alpha}(Q)h_{\alpha}(S) },\label{schottky}
\end{equation}
where $\theta[\alpha](\boldsymbol{u})$ is a $\theta$-function with
non-singular odd characteristics $[\alpha]$, $\mathrm{d}
\boldsymbol{v}$ is the vector of normalized holomorphic differentials,
and
\[
h_{\alpha}(Q)^2=\sum_{k=1}^g \frac{\partial}{\partial z_k}
\theta[\alpha](\boldsymbol{0})\,\mathrm{d} v_k(Q)
\]
The bi-differential $\omega(Q,S)$ is then given by \cite{fa73}.
\begin{equation}
\omega(Q,S)=\mathrm{d}_Q\mathrm{d}_S \,\mathrm{ln}\, E(Q,S).
\end{equation}

We emphasize that in this paper, we will instead rely on alternative
``algebraic'' constructions of the differential $\omega(Q,S)$.  By
following classical works such as \cite{kl86,kl88}, together with
results documented in \cite{ba97} we express the differential
$\omega(Q,S)$ in the form
\begin{equation}
  \omega(Q,S)=\frac{\mathcal{F}(Q,S)}{f_y(Q) f_w(S)(x-z)^2}
  \mathrm{d}x\mathrm{d}z+
  2\,\mathrm{d}\boldsymbol{u}(Q)^T\varkappa \mathrm{d}
  \boldsymbol{u}(S),
\end{equation}
where $Q=(x,y)$, $S=(z,w)$, and the function
$\mathcal{F}(Q,S)=\mathcal{F}((x,y),(z,w))$ is a polynomial of its
arguments, with coefficients depending on the moduli of the curve
$X_g$.  Finally, $\varkappa$ is a symmetric matrix $\varkappa^T=\varkappa$
that is chosen to provide a normalization of $\omega(Q,S)$; it is
expressible in terms of the first and second period matrices
$\varkappa= \omega^{-1}\eta$. We will refer to the term of
$\omega(Q,S)$ including the polynomial $\mathcal{F}(Q,S)$ as its
algebraic part,
\begin{align}
\omega(Q,S)=\omega^{\rm{alg}}(Q,S)+
2\, \mathrm{d}\boldsymbol{u}(Q)^T\varkappa \mathrm{d}\boldsymbol{u}(S)\label{alg}
\end{align}
In the vicinity of a base point $P$, where points $Q$ and $S$ are represented
by local
coordinates $\xi(Q)$ and $\xi(S)$ respectively, the holomorphic part
of $\omega^{\rm{alg}}(Q,S)$ is expanded in the series
\begin{align}\begin{split}
  &\left.\frac{\mathcal{F}(Q,S)\mathrm{d}z \mathrm{d}x }{f_y(Q) f_w(S)
      (x-z)^2} \right|_{ x=x(Q),z=x(S)} - \frac{\mathrm{d} \xi(Q)
    \mathrm{d} \xi(S)}{(\xi(Q)-\xi(S))^2 }
  \\&=\sum_{k,l=0}^{\infty}\omega_{k,l}^{\mathrm{alg}}(P) \xi(Q)^k
  \xi(S)^l \mathrm{d} \xi(Q) \mathrm{d} \xi(S).\end{split} \label{omegaalg}
\end{align}

An algorithm to construct the polynomial $\mathcal{F}(Q,S)$ is known,
see e.g.\ \cite{ba97} and therefore functions such as
$\omega_{k,l}^{\mathrm{alg}}(P)$ are needed for our construction are
considered as known.  In all that follows we will take the fixed base point
$P$ to be $(\infty,\infty)$.  We shall present below some explicit
expressions for $\mathcal{F}$ as well as the first few terms of the
expansions $\omega_{k,l}^{\rm{alg}}(P)$ in the simplest cases:

{\sf Example I}:

\begin{equation}\label{bi-diff}
\omega^{\rm{alg}}(Q,S)=\frac{F(x,z)+2yw}{4(x-z)^2}
\frac{\mathrm{d}x}{y}\frac{\mathrm{d}z}{w},\quad Q=(x,y), S=(z,w),
\end{equation}
where
\begin{equation} \label{hypbidiff}
F(x,z)= 4 x^2 z^2 (x+z) + 2 \alpha_4 x^2 z^2 +\alpha_3
x z(x+z) + 2 \alpha_2 xz +\alpha_1(x+z) +2 \alpha_0.
\end{equation}
The polynomial $F(x,z)$, with the properties $F(x,z)=F(z,x)$ and
$F(x,x)=2y^2$, is sometimes called the  ``Kleinian polar''.

Expanding this about the base point gives:
\begin{align*}
\omega_{0,0}^{\rm{alg}}&=-\frac{\alpha_4}{8},\\
\omega_{0,1}^{\rm{alg}}&=\omega_{1,0}^{\rm{alg}}=0,\\
\omega_{0,2}^{\rm{alg}}&=\omega_{2,0}^{\rm{alg}}=
-\frac{16\alpha_3-3\alpha_4^2}{128},\qquad
\omega_{1,1}^{\rm{alg}}=0,\\
&\ldots\text{etc.}
\end{align*}
\begin{remark}
For hyperelliptic curves, the coefficients $\omega_{i,j}^{\rm{alg}}$
vanish if either of $i$ or $j$ is odd.
\end{remark}

{\sf Example II}:

 \begin{equation}
   \omega^{\rm{alg}}((x,y),(z,w))=\frac{\mathcal{F}((x,y),(z,w))\,dxdz} {(x-z)^2
     f_y(x,y)f_w(z,w)}\label{omegatrig}
\end{equation}
with the polynomial $\mathcal{F}((x,y);(z,w))$ given by the formula
\begin{align}\begin{split}
\mathcal{F}\big(&(x,y),(z,w)\big)=w^2 y^2\\
&+w\left(w \left[\frac{f(x,y)}{y}\right]_y +T(x,z)\right)
 +y\left(y \left[\frac{f(z,w)}{w}\right]_w +T(z,x)\right)\end{split}
 \label{omegatrig1}
\end{align}
and
\begin{align}
    T(x,z)&=3\mu_{12}+ ( z+2x ) \mu_9+x ( x+2\,z )\mu_6\\
    &+3\mu_3 x^2 z + x^2 z^{2}+2\,x^3z.\label{omegatrig2}
\end{align}

Expanding this about the base point gives:
\begin{eqnarray*}
  \omega_{0,0}^{\rm{alg}} =&0, &\\
  \omega_{0,1}^{\rm{alg}} =&\omega_{1,0}^{\rm{alg}}&=-\frac{2}{3}\mu_3,\\
  \omega_{0,4}^{\rm{alg}}=&\omega_{4,0}^{\rm{alg}} & =-\frac{2}{3}\mu_6
  +\frac{5}{9}\mu_3^2,\\
  \omega_{1,3}^{\rm{alg}} =&\omega_{3,1}^{\rm{alg}} & =-\frac{2}{3}\mu_6
  +\frac{4}{9}\mu_3^2,\\
  \omega_{2,2}^{\rm{alg}} =&0, &\\
  \ldots\text{etc.}
\end{eqnarray*}

\begin{remark} $\omega_{i,j}^{\rm{alg}}=0$ unless $i+j\equiv 1\,
  \mathrm{mod}\, 3$. This is a consequence of the cyclic symmetry of
  the curve. \end{remark}

\section{Algebro-geometric $\theta$, $\sigma$ and $\tau$-functions}

Let $\boldsymbol{v}\in \mathrm{Jac}(X_g) $ and $\theta(\boldsymbol{v})$
be a canonical $\theta$-function, that is, a $\theta$-function with zero
characteristics:
\[
\theta(\boldsymbol{v})=\sum_{\boldsymbol{m}\in \mathbb{Z}^g} \mathrm{exp}
\left\{ \imath\pi \boldsymbol{m}^T\tau \boldsymbol{m} +2\imath \pi
\boldsymbol{v}^T \boldsymbol{m} \right\}.
\]
The starting point of this paper is the transition from $\theta$ to
$\sigma$-functions. For any point $\boldsymbol{u}\in \mathrm{Jac}(X)$
we define:
\begin{equation}
\sigma(\boldsymbol{u})=C(\tau)\theta(\mathcal{A}^{-1} \boldsymbol{u})
\mathrm{exp}\left\{  \frac12 \boldsymbol{u}^T \varkappa \boldsymbol{u}
 \right\},
\end{equation}
where $\theta(\boldsymbol{v})$ is the canonical $\theta$-function, and
$C(\tau)$ is a certain modular constant that we do not need for the
results that follow.  We note that this $\sigma$-function differs from
the ``fundamental $\sigma$ function'' of the publications mentioned in
the introduction by the absence of a shift of the $\theta$-argument by
the vector of Riemann constants. Thus $\sigma(\boldsymbol{0})\neq 0$.
The sigma-function inherits a quasi-periodicity property from the
corresponding $\theta$-function.

We introduce the Kleinian multi-variable $\zeta$ and $\wp$-functions as above
(\ref{zeta}, \ref{Kleinwp}). These functions are suitable coordinates
to describe Abelian functions and the KP-type hierarchies of
differential relations between them.

In higher genera, the relations between Abelian functions become
somewhat lengthy. These relations can often be summarized concisely by
developing a matrix formulation of the theory, as was done \cite{bel97b}
in the hyperelliptic case.

The Sato-Fay algebro-geometric $\tau$-function of the genus $g$ curve
$X_g$ of arguments $\boldsymbol{t}=(t_1,\ldots,t_g,t_{g+1},\ldots )^T$,
$\boldsymbol{u}\in \mathrm{Jac}(X)$, $P\in X_g$ is defined as
\begin{align}
\tau(\boldsymbol{t};\boldsymbol{u},P)=\theta\left( \sum_{k=1}^{\infty}
\boldsymbol{U}_k(P) t_k+\mathcal{A}^{-1}\boldsymbol{u} \right) \mathrm{exp}
\left\{ \frac12\sum_{m,n\geq 1} \omega_{mn}(P) t_mt_n \right\}.\label{tau}
\end{align}
Here $\mathcal{A}$ is the matrix of periods of canonical holomorphic
differentials and the winding vectors defined in (\ref{winvec}),
$\boldsymbol{U}_k(P)$, appear in the expansion of the normalized
holomorphic integral $\boldsymbol{v}$ in the vicinity of the given
point $P\in X$,
\[
\int\limits_{P_0}^Q \mathrm{d}\boldsymbol{v}(Q')
=\int_{P_0}^P\mathrm{d}\boldsymbol{v}(Q')+
\sum_{k=1}^{\infty} \boldsymbol{U}_k(P) \xi(Q)^k
\]
with $\xi(Q)$ being the local coordinate of the point $Q$ in the
vicinity of the given point $P$, so that $\xi(P)=0$. The quantities
$\omega_{mn}(P)$ define the holomorphic part of the expansion of the
fundamental second kind differential $\omega(Q,S)$ near the point $P$
according to (\ref{omegapole}).

We restrict ourselves to the case of algebraic curves with a branch
point at infinity, and take $P=(\infty,\infty)$ to be the base point
where we expand all our functions. The winding vectors are in this case
\begin{equation}
\boldsymbol{U}_k(\infty)\equiv\boldsymbol{U}_k=\mathcal{A}^{-1}
\boldsymbol{R}_k,\quad k=1,\ldots,g, \label{winvec}
\end{equation}
where $\boldsymbol{R}_1,\boldsymbol{R}_2,\ldots$ are residues of
canonical holomorphic integrals multiplied by differentials of the
second kind with poles of order $k$ at infinity, giving:
\[
\boldsymbol{R}_k= \frac{1}{k}\left. \frac{\mathrm{d}^{k-1}}{ \mathrm{d}
\xi(Q)^{k-1}}\mathrm{d} \boldsymbol{u}(Q) \frac{1}{\mathrm{d}\xi(Q) }
\right|_{Q=\infty}.
\]
Using the above definitions, we can see that the algebro-geometric
$\tau$-function is given by the formula,
equivalent to (\ref{tausigma}):
\begin{align}
\frac{ \tau(\boldsymbol{t};\boldsymbol{u})}
{\tau(\boldsymbol{0};\boldsymbol{u})}=
\frac{ \sigma\left(\sum_{k=1}^{\infty} \boldsymbol{R}_k t_k+
\boldsymbol{u}\right)}
{\sigma(\boldsymbol{u})}\,\mathrm{exp} \left\{ \frac12
\sum_{k,l=0}^{\infty} \omega_{k,l}^{\mathrm{alg}}t_kt_l  \right\}.
\label{tausigma2}\end{align}
In (\ref{tausigma2}), $\omega_{k,l}^{\mathrm{alg}}$ is the
\emph{algebraic part} of the holomorphic part $\omega_{k,l}$ of the
expansion of the bi-differential as defined in (\ref{omegaalg}).  One
can see that the non-algebraic part, i.e.\ the normalizing bi-linear
form, is absorbed into the $\sigma$-function.

Once we have an expression involving derivatives of the sigma
function, we need to convert this to an expression involving
derivatives of the $\wp$-function.  To do this we start with the
definition of the $\zeta$ function (\ref{zeta}), which we write in the form
\[
\sigma_i(\boldsymbol{u})=\zeta_i(\boldsymbol{u})\sigma(\boldsymbol{u}),
\quad i=1,\ldots,g
\]
then repeated differentiation gives us a ladder of relations which
enable us to recursively express any derivative of $\sigma$ in terms
of $\wp_{ij\dots k}$, $\zeta_i$ and $\sigma$, all evaluated at
$\mathbf{u}$.
\begin{align}
  \sigma_{ij}(\boldsymbol{u})&=\sigma_j(\boldsymbol{u})\zeta_i(\boldsymbol{u})-
  \sigma(\boldsymbol{u}) \wp_{ij}(\boldsymbol{u}) \label{sigmaladder}\\
  \sigma_{ijk}(\boldsymbol{u})& =\sigma_{jk}(\boldsymbol{u})\zeta_i(\boldsymbol{u})
  -\sigma(\boldsymbol{u}) \nonumber
  \wp_{ijk}(\boldsymbol{u})-\sigma_k(\boldsymbol{u})
  \wp_{ij}(\boldsymbol{u}) -\sigma_j(\boldsymbol{u})
  \wp_{ik}(\boldsymbol{u}), \\\text{etc} \dots \nonumber
\end{align}

In the following sections, we will look at different methods for constructing
such relations between the derivatives of $\sigma$, and hence between the
Abelian functions associated with the curve.

\section{The ``classical'' method}

The starting point for this approach is the Klein formula, which compares
two different expressions for the fundamental bi-differential:
\begin{theorem}
  Let the canonical holomorphic differentials of the curve $f(x,y)=0$ be
  represented in the form
\[
\mathrm{d} u_i(x,y)=\frac{\mathcal{U}_i(x,y)}{ f_y(x,y) } \mathrm{d}x,\quad i=1,\ldots,g,
\]
where $\mathcal{U}_k(x,y)$ are monomials of their variables.
Then the following $g$ relations hold
\begin{align}
\label{klein}\begin{split}
  &\sum_{i,j=1}^g \wp_{ij}\left( \int_{P_0}^{(x,y)}\mathrm{d}
    \boldsymbol{u}-\sum_{k=1}^g \int_{P_0}^{(x_k,y_k)}\mathrm{d}
    \boldsymbol{u}+\boldsymbol{K}_{P_0} \right)
  \mathcal{U}_i(x,y)\mathcal{U}_j(x_k,y_k)\\&=
  \frac{\mathcal{F}(x,y;x_k,y_k)}{(x-x_k)^2},\quad k=1,\ldots,g.
\end{split}
\end{align}
with polynomial $\mathcal{F}(P,Q)=\mathcal{F}(x,y;z,w)$ defined in
(\ref{alg}), an arbitrary base point of the Abel map $P_0$, and the
corresponding vector of Riemann constants, $\boldsymbol{K}_{P_0}$.

\end{theorem}

As an application of this theorem consider the case of a curve with
branch point at $P_0=(\infty,\infty)$. Let us tend $P=(x,y)\rightarrow
(\infty,\infty)$. Both, right and left had sides of (\ref{klein}) have
poles at infinity. Equating principal part of poles, we get a set of
relations between $P_k=(x_k,y_k)$ and the multi-indexed $\wp$-symbols
which can be interpreted as differential equations for the
$\wp_{k,l}$-functions. In this way, the solution of the Jacobi
inversion problem can also be derived in terms of $\wp$-functions.

\subsection{ Example: Hyperelliptic curve of genus two}
The $\sigma$-functional realization of hyperelliptic functions of a
genus two curve has already been discussed in many places, see
e.g. \cite{ba97,bel97b}. But we shall briefly describe here
the principal points of the construction to convey the structure of
the theory that we wish to develop for higher genera curves.

We consider the genus two hyperelliptic curve of {\sf example I}.
The algebraic part of the fundamental bi-differential is given as
above in eqn.  (\ref{bi-diff}) \cite{ba97}.

Following the procedure in those papers, we equate principal parts of the
highest (second order) pole to obtain
\begin{align}
\wp_{12}+x \wp_{11}-x^2=0,\label{JIP2}
\end{align}
which is the $x$-part of the Jacobi inversion problem for this curve.
Then equating principal parts of the next, lower pole and using
(\ref{JIP2}) leads to the relation
 \begin{align}
y_k=-\wp_{112}-x_k \wp_{111},\label{JIP2A}
\end{align}
which completes the solution of the Jacobi inversion problem.

We can use these relations to eliminate $y$ and quadratic and higher
terms in $x$, and we get from the next term
\begin{align*}
  &(\tfrac12 \wp_{1111} -3 \wp_{11}^2-\tfrac12 \alpha_4 \wp_{11} - 2
  \wp_{12} - \tfrac14 \alpha_3 ) x_k\\ &\qquad+ \tfrac12 \wp_{1112} - \tfrac12
  \alpha_4 \wp_{12} - 3 \wp_{11} \wp_{12} + \wp_{22}=0.
\end{align*}
This equation holds for both $x_k$, so the coefficients of different
powers of $x_k$ must be identically zero, and solving for
$\wp_{1111}$ and $\wp_{1112}$ we find
\begin{align}
 \wp_{1111} &= 6 \wp_{11}^2 + \alpha_4 \wp_{11} + 4 \wp_{12} +
  \tfrac12 \alpha_3,\label{R1111}\\
 \wp_{1112} &= 6\wp_{11} \wp_{12} +\alpha_4 \wp_{12}  -2\wp_{22},\label{R1112}
\end{align}
the first two 4-index relations in the genus 2 case.  A relation
for 3-index functions $\wp_{ijk}$ can be obtained by cross-derivation of
(\ref{R1111}) and (\ref{R1112}):
\[
\wp_{122}+\wp_{11}\wp_{112}-\wp_{12}\wp_{111}=0.
\]

Other equations are found from higher order terms.  At every stage we
need to substitute for higher derivative terms (such as $\wp_{11111}$,
for example) by using derivatives of previously derived relations. In
addition, multiplication by a 3-index $\wp_{ijk}$ is sometimes useful,
followed by substitution of known relations which are quadratic in the
$\wp_{ijk}$. The first such relation is found to be
$$
\mathrm{Jac}_{6}:\quad \wp_{111}^2  = 4 \wp_{11}^3+\alpha_{3}
\wp_{11}+\alpha_{4} \wp_{11}^2+4 \wp_{12} \wp_{11}+\alpha_{2}+4
\wp_{22}.
$$
This is one of the relations that describe the Jacobi variety of the
curve as algebraic variety.  These relations are well documented in
\cite{ba07}, see also \cite{bel97b}. Below we show an alternative way
to derive these relations.


We remark briefly that one variation on the classical method is to
invoke the use of the sigma function expansion in the $u_i$ variables.
The first few terms in the sigma expansion are derived from the lowest
order expansion terms of (\ref{klein}).  These can then be used in a
bootstrap fashion to derive higher order PDEs for the $\wp$ functions.
This approach was essential in the genus six (4,5) case considered in
\cite{ee09}, which is discussed in more detail in that paper.


\section{Derivation of integrable hierarchies via Pl{\"u}cker
  relations}
The key to this approach is Sato's formula, Theorem (5.1) below.  This
gives an expansion of a ratio of $\tau$-functions,
${\tau(\boldsymbol{t};\boldsymbol{u})}/
{\tau(\boldsymbol{0};\boldsymbol{u})}$, as a series of rational
expressions in the $\tau$-function and its derivatives
at$\boldsymbol{t}=0$. The $\boldsymbol{t}$-dependence is given in
terms of Schur polynomials $s_\lambda(\boldsymbol{t})$ in the times
$t_i$. The coefficients of these polynomials are determinants of
differential expressions of $\tau$, which are Pl\"ucker coordinates on
a Grassmannian. Such coordinates satisfy the Pl\"ucker relations -
each partition $\lambda$ can be expanded in hooks, and the
corresponding Pl\"ucker coordinates are expressible, analogously to
Giambelli's formula, as determinants of single hook partitions. These
relations give the differential relations for the Abelian functions
which we seek.

For any partition $\lambda: \alpha_1\geq \alpha_2\geq \ldots \geq
\alpha_n$ of $|\lambda|=\sum_{i=1}^n\alpha_i$, the Schur polynomial of
$n$ variables $x_1,\ldots,x_n$ is defined by
\[
s_{\lambda}(\boldsymbol{x})=\mathrm{det}\left( p_{\alpha_i-i+j}(\boldsymbol{x})
\right)_{i,j=1,\ldots,n},
\]
where the elementary Schur functions $p_m(\boldsymbol{x})$ are
generated by the series,
\[
\sum_{m=0}^{\infty} p_m(\boldsymbol{x}) t^m =\mathrm{exp}
\left\{ \sum_{n=1}^{\infty} x_n t^n\right\}.
\]
The first few Schur polynomials are
\begin{align*}
  &s_{1}(\boldsymbol{x})=x_1, \\
  & s_{2}(\boldsymbol{x}) =x_2+\tfrac12
  x_1^2,\quad
  s_{1,1}(\boldsymbol{x})=-x_2+\tfrac12 x_1^2,\\
  &s_{3}(\boldsymbol{x})=x_3+x_1x_2+\tfrac16 x_1^3, \quad
  s_{2,1}(\boldsymbol{x})=-x_3+\tfrac13 x_1^3,\quad
  &s_{1,1,1}(\boldsymbol{x})=x_3-x_1x_2+\tfrac16 x_1^3\\
  &s_{4}(\boldsymbol{x})=x_4+x_1x_3+\tfrac12 x_2^2+\tfrac12 x_1^2x_2
  +\tfrac{1}{24}x_1^4, \quad \text{etc.}
\end{align*}

The Cauchy-Littlewood formula
\[
\mathrm{exp}\left\{ \sum_{n=1}^{\infty} n x_n y_n \right\}
=\sum_{\lambda}s_{\lambda}(\boldsymbol{x})
s_{\lambda}(\boldsymbol{y}),
  \]
  where $s_{\lambda}(\boldsymbol{x})$ is the Schur function of the partition
  $\lambda:\alpha_1\geq \alpha_2\geq \ldots \geq \alpha_n$, leads to
  the Taylor expansion in the form
\begin{align*}
  f(\boldsymbol{x})=\left.\mathrm{exp}\left\{
      \sum_{n=1}^{\infty}x_n\frac{\partial}{\partial y_n} \right\}
    f(\boldsymbol{y}) \right|_{\boldsymbol{y}=0}
  =\left.\sum_{\lambda}s_{\lambda}(\boldsymbol{x})
    s_{\lambda}\left(\frac{1}{n} \frac{\partial}{\partial y_n} \right)
    f(\boldsymbol{y})\right|_{\boldsymbol{y}=0}
\end{align*}

We now introduce the Frobenius notation for partitions.  Single hook
partitions $(\alpha,\beta)$ are denoted $(\alpha+1,1^{\beta})$. All
partitions can be decomposed into finitely many hooks:
\[
(\lambda)=(\alpha_1,\ldots,\alpha_r\vert \beta_1,\ldots,\beta_r ).
\]
The total number of hooks $r$ is called the {\em rank} of the partition
$\lambda$, \cite{sag01}.
In particular, a partition which decomposes into two hooks
(rank 2) is written as:
\[
(n,m, 2^k, 1^l) = (n-1, m-2 \vert k+l+1,k),
\]
where $n>m>1$, $k\ge 0$ and $l\ge 0$.  Most of the formulae below are
derived using rank 2 partitions.

Giambelli's formula now shows how to expand a Schur function $s_{\lambda}$
in hooks,
\[
s_{\lambda}(\boldsymbol{x})=\det \left(
s_{(\alpha_i,\beta_j)}(\boldsymbol{x}) \right)_{1\leq i,j\leq r},
\]

\begin{theorem}[\bf Sato formula ] Let
  $\tau(\boldsymbol{t};\boldsymbol{u})$ be any function of vector
  arguments
\[
\boldsymbol{U}_1t_1,\boldsymbol{U}_2t_2,\ldots,\quad
\boldsymbol{t}=(t_1,t_2,\ldots)\in\mathbb{C}^{\infty},
\]
where $\boldsymbol{U}_1,\boldsymbol{U}_2\ldots $ is an infinite set of
constant complex vectors from $\mathbb{C}^{g}$ and
$\boldsymbol{u}\in\mathbb{C}^g$ is a parameter. Suppose that
$\tau(\boldsymbol{0};\boldsymbol{u})\neq 0$. Then for any partition
$(\lambda)=(\alpha_1,\ldots,\alpha_r\vert \beta_1,\ldots,\beta_r) $
and any $\boldsymbol{u}\in\mathbb{C}^g$
\begin{equation}
\frac{\tau(\boldsymbol{t};\boldsymbol{u})}
{\tau(\boldsymbol{0};\boldsymbol{u})}=
\sum_{\lambda} s_{\lambda}(\boldsymbol{t})\det \left( (-1)^{\beta_j+1}
A_{(\alpha_i|\beta_j)} (\boldsymbol{u}) \right),
\end{equation}
where the $A_{(m|n)}(\boldsymbol{u})$ with  $m\geq 0$, $n\geq 0$ form a
linear basis of the Grassmannian:
\begin{align}\begin{split}
    A_{(m|n)}(\boldsymbol{u})&=-A_{(n|m)}(-\boldsymbol{u})\\
    &=\left.(-1)^{n+1}s_{m+1,1^n}
      (\partial_{\boldsymbol{t}})\tau(\boldsymbol{t};\boldsymbol{u})
    \right|_{\boldsymbol{t}=0}\tau(\boldsymbol{0};\boldsymbol{u})^{-1}\\
    &=\left.\sum_{\alpha=0}^m p_{n+\alpha+1}(-\partial
      _{\boldsymbol{t}}) p_{m-\alpha}(\partial_{\boldsymbol{t}})
      \tau(\boldsymbol{t}; \boldsymbol{u})\right|_{\boldsymbol{t}=0}
    \tau(\boldsymbol{0}; \boldsymbol{u})^{-1}
\end{split}\label{Amn}
\end{align}
and
\[
\partial _{\boldsymbol{t}} =\left( \frac{\partial}{\partial t_1},
  \frac12 \frac{\partial}{\partial t_2}, \frac13
  \frac{\partial}{\partial t_3},\ldots, \right).
\]
\end{theorem}

\begin{theorem}[\bf Pl\"ucker relations] For any partition
  $(\lambda)=(\alpha_1,\ldots,\alpha_r\vert \beta_1,\ldots,\beta_r )$
  and any $\boldsymbol{u}\in\mathbb{C}^g$, the $\tau$-function satisfies:
\begin{equation}
  \left.\tau(\boldsymbol{0};\boldsymbol{u})^{r-1}s_{\lambda}
    (\partial_{\boldsymbol{t}})\tau(\boldsymbol{t};\boldsymbol{u})
  \right|_{\boldsymbol{t}=0}=\left.\det \left( s_{(\alpha_{i}|\beta_j)}
      (\partial_{\boldsymbol{t}}) \tau(\boldsymbol{t};\boldsymbol{u})
    \right|_{\boldsymbol{t}=0}  \right).\label{pr}
\end{equation}
\end{theorem}
In particular, for any curve $X$, its Sato algebro-geometric $\tau$-function
has the expansion
\begin{align}
\begin{split}
  \frac{\tau(\boldsymbol{t};\boldsymbol{u})}
  {\tau(\boldsymbol{0};\boldsymbol{u})}& =1+ A_{(0|0)}(\boldsymbol{u})
  {s_1} (\boldsymbol{t})
  +A_{(1|0)}(\boldsymbol{u}){s_{2}}(\boldsymbol{t})\\
&+A_{(0|1)}(\boldsymbol{u})
  {s_{1,1}}(\boldsymbol{t})+\ldots\\
  &=1+ A_{(0|0)} t_1 + A_{(1|0)} (t_2 +\tfrac12 t_1^2)
  +A_{(0|1)}(\boldsymbol{u})
( -t_2+\tfrac12 t_1^2)+\ldots,
\end{split}
\end{align}
where $A_{(m|n)}(\boldsymbol{u})$ as defined in (\ref{Amn}) are the elements
of the
semi-infinite matrix $A$,
\[
A=\left(\begin{array}{cccc} A_{(0|0)}(\boldsymbol{u})
&A_{(0|1)}(\boldsymbol{u})&A_{(0|2)}(\boldsymbol{u})
&\ldots\ldots\\
    A_{(1|0)}(\boldsymbol{u})&A_{(1|1)}(\boldsymbol{u})
&A_{(1|2)}(\boldsymbol{u})&\ldots\ldots\\
    A_{(2|0)}(\boldsymbol{u})&A_{(2|1)}(\boldsymbol{u})
&A_{(2|2)}(\boldsymbol{u})&\ldots\ldots\\
    \vdots&\vdots&\vdots&\ldots\ldots
  \end{array} \right)  = (\boldsymbol{A}_0,
\boldsymbol{A}_1,\ldots\ldots)
\]
with infinite vectors $\boldsymbol{A}_i$, $i=1,2,\ldots$. To more
complicated partitions such as $(m_1,\ldots,m_k|n_1,\ldots,n_k)$ we
associate according to the Giambelli formula certain minors of $A$,
 \[
 A_{(m_1,\ldots,m_k|n_1,\ldots,n_k)}(\boldsymbol{u})=\mathrm{det}
\left( A_{(m_i|n_j  ) (\boldsymbol{u})} \right)_{i,j=1\ldots,k}.
 \]
 Note that  these $A$'s are not independent -- they satisfy the
 Pl\"ucker relations. These are a family of differential equations
 satisfied by $\tau$, that represent  a completely integrable hierarchy
 of KP-type.

 The definitions and relations given above are valid for any
 multivariate function $\tau(\boldsymbol{t};\boldsymbol{u})$. Below
 we consider functions $\tau$-functions constructed on the Jacobi
 varieties of algebraic curves.

To each symbol  we shall put in correspondence its weight
$$
\wp_{\underbrace{\scriptstyle 1,\ldots, 1}_{k_1},
  \underbrace{\scriptstyle 2,\ldots, 2}_{k_2}, \ldots,
  \underbrace{\scriptstyle g,\ldots, g}_{k_g}}\quad \Leftrightarrow
\quad \sum_{j=1}^gk_jw_j,
$$
where $w_i$ is the order of vanishing of the holomorphic integral
$\int\mathrm{d}u_i$ at infinity. In other words, if the curve has
a Weierstrass point at infinity, then $w_i=1<w_2<\ldots<w_g$ is the
Weierstrass gap sequence at infinity.

For a given weight $W$ consider all Young diagrams which
decompose into only two hooks, and write
corresponding relations between multi-index symbols
$\wp_{i_1,\ldots,i_n}$. The first non-trivial Young diagram
corresponds to the partition $\lambda=(2,2)$. There is only one
multi-index symbol of weight 4, that is $\wp_{1111}(\boldsymbol{u})$
and its corresponding Pl\"ucker relation for any curve is of the form
\begin{equation}
\wp_{1111}(\boldsymbol{u})=\text{ polynomial of even symbols} \;
\wp_{ij}(\boldsymbol{u})
\end{equation}
For higher weights a larger number of multi-index functions can be
constructed, and this number grows rapidly with increasing weight. To
find them in terms in the form of polynomials of two-index functions,
we shall write Pl\"ucker relations corresponding to independent
diagrams of the same weight, and solve the corresponding linear
systems.  The technique is best illustrated by examples.

We consider the genus two hyperelliptic curve, with the differentials and
the bi-differential chosen as in {\sf example I} above.
The simplest class of non-trivial Pl\"ucker relations is found from the class
of Young diagrams which may be decomposed into two hooks, that is, those
of the form $(2+m,2+n,2^k,1^l)$, with $m\ge n \ge 0, k \ge 0, l \ge 0$. In
Frobenius  `hook' notation these are $(m+1,n|k+l+1,l)$.
The corresponding Pl\"ucker relations read
\begin{equation}
A_{(m+1,n|k+l+1,l)} = \left|\begin{array}{cc} A_{(m+1,k+l+1)}& A_{(m+1,l)}\\
                                              A_{(n,k+l+1)}& A_{(n,l)}
\end{array}\right|.
\label{two-hook}
\end{equation}
These equations are all bilinear partial differential equations in the
$\tau$-function.  They may be expanded in terms of the Kleinian
$\wp_{ij}$ and $\zeta_i$ functions.

The first Young diagram leading to a non-trivial Pl\"ucker relation
corresponds to the partition $\lambda=(2,2)$, with Young diagram
$$ \yng(2,2).$$
Writing (\ref{pr}) in
this case we obtain after simplification
\begin{equation}
\mathrm{KdV}_{4}: \qquad \wp_{1111}(\boldsymbol{u})=
6\wp_{11}^2(\boldsymbol{u})+4\wp_{12}(\boldsymbol{u})
+\alpha_4\wp_{11}(\boldsymbol{u})+\tfrac12\alpha_3. \label{p1111}
\end{equation}
This is of course the same equation as (\ref{R1111}).
The weight of the diagram is 4 in this case, which is the same as the
weight of the equation, defined as the weighted sum of indices in
which the index ``1'' has weight 1 and index ``2'' has weight 3.  We
will use this weight correspondence in other cases too, and will
denote the weight of the object by subscript $i+3 j$ where $i$ and $j$
are respectively the numbers of 1's and 2's in the multi-index
relation.
The relation (\ref{p1111}) is the analogue of Weierstrass' equation for $\wp''$
in the genus 1 case.

The next group of diagrams are of weight 5, and correspond to the
partitions $\lambda=(3,2)$ and $\lambda=(2,2,1)$, with their
transposes, which give the same equations, in the hyperelliptic case.
Both these Pl\"ucker relations lead to the equation
\[
\left( \zeta_1(\boldsymbol{u}) + \frac{\partial}{\partial u_1}\right)
\mathrm{KdV}_{4}=0
\]
i.e. the Young diagram of weight 5 gives no new equations and we conclude
that the following correspondence is valid
\[
\left\{ \begin{array}{c}  \lambda=(2,2)\\
    \lambda=(3,3),\quad \lambda=(2,2,1) \end{array}\right\}
\Longleftrightarrow \mathrm{KdV}_{4}
\]
At weight 6 we have three independent Young diagrams with $(2,2)$ centres
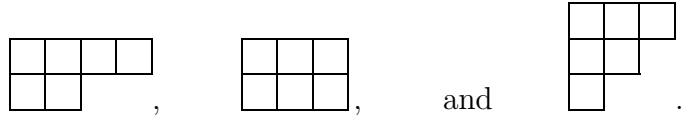
\[
\yng(4,2),\hskip1cm \yng(3,3)\,,\hskip1cm \text{and} \hskip1cm \yng(3,2,1).
\]
The other two diagrams of weight 6
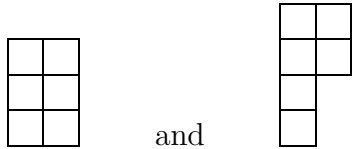
\[
\yng(2,2,2) \hskip1cm \text{and} \hskip1cm\yng(2,2,1,1)
\]
again give the same results  as their transposes.
\begin {remark}The Pl\"ucker relations associated
with any Young diagram and its transpose are always the same for a
hyperelliptic curve; however for general curves, this is no longer true.
\end{remark}

These three independent diagrams give an overdetermined system of three
equations.  After substituting for higher derivatives of the known
relation (\ref{p1111}), we can solve for the two unknowns $\wp_{1112}$
and $\wp_{111}^2$ to get
\begin{align}
  \mathrm{KdV}_{6}:\quad \wp_{1112} &= 6 \wp_{12} \wp_{11}-2 \wp_{22}
  + \alpha_{4} \wp_{12}\label{p1112}\\
  \mathrm{Jac}_{6}:\quad \wp_{111}^2 & = 4 \wp_{11}^3+\alpha_{3}
  \wp_{11}+\alpha_{4} \wp_{11}^2+4 \wp_{12} \wp_{11}+\alpha_{2}+4
  \wp_{22}\label{p111s}
\end{align}
The relation (\ref{p1112}) is another generalization of the equation
for $\wp''$ in the genus one case, and is identical to (\ref{R1112}) found
using the classical method. The relation (\ref{p111s})
corresponds to Weierstrass' equation for $(\wp')^2$ in the genus one
case.  In what follows we will always assume that higher derivatives of the
known relations at a lower weight have been eliminated.

At weight 7 we have 5 independent Young diagrams with (2,2) centres.
The resulting overdetermined system of equations, of weight 7, contain
$\zeta_1$ multiplied by linear combinations of the two weight 6
relations (\ref{p111s},\ref{p1112}).  In addition we get a
relation at weight 7, namely
\[
\wp_{12} \wp_{111}-\wp_{122}-\wp_{112} \wp_{11}=0, \label{L5}
\]
as noted above.  We refer to relations of this type, linear in the
three-index $\wp_{ijk}$, as quasilinear.  These have no counterpart in
the genus 1 theory.  They can also be derived by cross-differentiation
between the two relations (\ref{p1111},\ref{p1112}), since
\[
\frac{\partial}{\partial u_2} \wp_{1111} = \frac{\partial}{\partial
  u_1} \wp_{1112}.
\]
Equations $KdV_4$ and $KdV_6$ describe the genus 2 solutions of the
KdV hierarchy associated with the given curve. To describe the Jacobi
variety and the Kummer variety we must consider diagrams of higher
weights.

For the diagram of weight 8 we find a set of eight overdetermined
equations with solution given by
\begin{eqnarray*}
\mathrm{ Jac}_8:\;\wp_{112}^2&=&\alpha_0-4\wp_{22}\wp_{12}+\alpha_4\wp_{12}^2
+4\wp_{11}\wp_{12}^2\\
\mathrm{KdV}_{8}:\;\wp_{1122}&=&2\wp_{11}\wp_{22}+4\wp_{12}^2+
\tfrac12\alpha_3\wp_{12},
\end{eqnarray*}
And at weight 9 the only new relation is the quasilinear relation
\[
8\wp_{122}\wp_{11}-4\wp_{12}\wp_{112}-4\wp_{222}+2\alpha_{4}\wp_{122}
-\alpha_{3}\wp_{112}-4\wp_{111}\wp_{22} = 0
\]
At weight 10 we have 18 overdetermined system of equations, which can
be solved for the 3 functions of weight 10, $\wp_{1222}, \wp_{112}^2,
\wp_{111}\wp_{122} $, giving
\begin{align*}
  \mathrm{KdV}_{10}:\;\wp_{1222}&=6\wp_{12}\wp_{22}+\alpha_2\wp_{12}
  -\tfrac12\alpha_1\wp_{11}-\alpha_0\\
  \mathrm{ Jac}_{10}^{(1)}:\;\wp_{111}\wp_{122}&=
  -\tfrac12\alpha_1\wp_{11} +2\wp_{22}\wp_{11}^2
  +2\wp_{11}\wp_{12}^2+\alpha_2\wp_{12}+4\wp_{22}\wp_{12}
  +\tfrac12\alpha_3\wp_{12}\wp_{11}\\
  \mathrm{ Jac}_{10}^{(2)}:\;\wp_{112}^2&=\alpha_0-4\wp_{22}\wp_{12}
  +\alpha_4\wp_{12}^2  +4\wp_{11}\wp_{12}^2
\end{align*}
The equations
$\mathrm{Jac}_8,\mathrm{Jac}_{10}^{(1)},\mathrm{Jac}_{10}^{(2)}$
represent an embedding of the Jacobi variety as a 3-dimensional
algebraic variety into the complex space $\mathbb{C}^5$ whose
coordinates are $\wp_{11},\wp_{12},\wp_{22},\wp_{111}, \wp_{112}$

At weight 12 we get the final 4-index relation for $\wp_{2222}$ and
two quadratic 3-index relations for $\wp_{112}\wp_{122}$ and
$\wp_{111}\wp_{222}$. We can continue in this manner at weight 14 to
get more quadratic 3-index relations \cite{bel97b}.  At the odd
weights 11, 13, 15, we get quasilinear relations which can also be
found by cross-differentiation.  As a practical point we note that the
equations we derive can often contain ideals generated by the lower
weight relations, and some work is required to identify genuinely new
relations.

At weight 16 we have 117 independent relations giving an
overdetermined system of equations (we have checked only a selection
of these). At this weight a new feature occurs.  As well as the
equations expressing the quadratic 3-index term $\wp_{122}\wp_{222}$
in terms of cubics in the $\wp_{ij}$, we have terms which are quartic
in the $\wp_{ij}$.  We can pick one of these quartic terms, say
$\wp_{12}^4$, and solve for this and for $\wp_{122}\wp_{222}$ to give
us two relations.  The relation involving a quartic in the $\wp_{ij}$
is just the Kummer variety of the curve. This is the quotient of the
Jacobi variety, $\mathrm{Kum}(X) = \mathrm{Jac}(X) /
(\boldsymbol{u}\rightarrow - \boldsymbol{u}$).  In the case $g=2$, the
Kummer variety is a surface in $\mathbb{C}^3$ which is given
analytically by a quartic equation.  This relation can can also be
found from the identity
\begin{equation}
(\wp_{111}^2)(\wp_{112}^2)-(\wp_{111}\wp_{112})^2=0. \label{kum25}
\end{equation}
The same quartic also appears, multiplied by various factors, at
higher weights.

\subsection{Example: Trigonal curve of genus three}

As before we consider {\sf Example II}, (2.3) whose holomorphic
differentials are given in (2.x)  \cite{eel00,bel00,eemop07}.  Here we
generally follow the notation of \cite{eemop07}). The fundamental
second kind differential is (\ref{omegatrig}), (\ref{omegatrig1}),
(\ref{omegatrig2}).

The first Young diagram leading to
a non-trivial Pl\"ucker relation, as in the genus 2 case, corresponds
to the partition $\lambda=(2,2)$.
Writing (\ref{pr}) in this case we obtain, after simplification,
\[
\wp_{1111}= 6\, \wp_{11}^{2}-3\,\wp_{22}.
\]
The weight of the diagram is again 4 in this case, which is the same
as the weight of the equation, defined as the weighted sum of indices
in which the index ``1'' has weight 1, index ``2'' has weight 2, and
index ``3'' has weight 5.  We will use this weight correspondence in
other cases too, and will denote the weight of the object by subscript
$i+2 j+5k$ where $i,j$ and $k$ are respectively the numbers of 1's,
2's and 3's in the multi-index relation. 

In the trigonal case we no longer longer have the symmetry about the
diagonal of the diagram that we have in the genus 2 case, but we can
restrict ourselves  by taking the symmetric or antisymmetric
combination of the two diagrams related by transposition.  In the
weight 5 case we have the antisymmetric combination
\[
\yng(3,2)\quad - \quad \yng(2,2,1),
\]
which gives the weight 5 trigonal PDE
\[
\wp_{1112}=6\,\wp_{11}\wp_{12}+3\,\mu_{3}\wp_{11}.
\]
For the symmetric case we get a derivative of the weight 4 equation,
plus $\zeta_1$ multiplied by the same equation.  With even (odd)
weights, the symmetric (antisymmetric) combinations give the 4-index
$\wp_{ijkl}$ relations.

At weight 6 we have the three symmetric diagrams
\begin{align*}
  \yng(4,2)\quad + \quad\yng(2,2,1,1), \quad \yng(3,3)\quad
  +\quad\yng(2,2,2)\,, \quad \yng(3,2,1),
\end{align*}
which give a set of three overdetermined equations with the unique
solution
\begin{align*}
  \wp_{111}^{2} & =4 \wp_{11} ^{3}+
  \wp_{12}^{2}+4\,\wp_{13}-4\,\wp_{11}\wp_{22},\\
  \wp_{1122}& =4\,\wp_{13}+2\,\mu_{6}+4\wp_{12}^{2}+
  2\,\wp_{11}\wp_{22}+3\,\mu_{3}\wp_{12}.
\end{align*}
Continuing in this way we recover the strictly trigonal versions of
the full set of equations given in \cite{eemop07}.

\section{Derivation of the integrable hierarchies via residues}
This method of deriving integrable equations on the Jacobians of
algebraic curves is based on the following relation due to Sato,
\cite{sato80}, \cite{sato81}, as applied to the algebro-geometric
context by Fay, \cite{fay83}.

The Hirota bilinear equations
\begin{eqnarray}
  &&\mathrm{Res}_{\xi=0} \frac{1}{\xi^2}
  \left[\mathrm{exp}\left\{\sum_{n=1}^{\infty} t_n\zeta^{-n} \right\}
    \exp\left\{- \sum_{n=1}^{\infty} \frac{\xi^n}{n}
      \frac{\partial }{\partial t_n} \right\}\tau(\boldsymbol{t};\boldsymbol{u})
  \right.\nonumber\\&& \left.
    \qquad\qquad\quad
    \mathrm{exp}\left\{\sum_{n=1}^{\infty} t_n'\xi^{-n} \right\}
    \exp\left\{- \sum_{n=1}^{\infty} \frac{\xi^n}{n}
      \frac{\partial }{\partial t_n'} \right\}\tau(\boldsymbol{t}';
    -\boldsymbol{u})\right]=0,\label{bilinearresidue}
\end{eqnarray}
where $\xi=\xi(Q)$ is the local coordinate of the point $Q$ near $P$, $\xi(P)=0$ defining the KP hierarchy are equivalent to the following Bilinear Identity

Denote $\Omega(\boldsymbol{f})$ the differential
\begin{equation}
\Omega(\boldsymbol{t},P)=\sum_{k=1}^{\infty} k t_k\Omega_{k}(P)
\end{equation}
where $\Omega_k(P)$ are normalized second kind differentials with
$k$-th order poles at infinity, so that $\Omega(\boldsymbol{t},P)$ has
an essential singularity
\[
\Omega(\boldsymbol{t},P) \simeq \left( \sum_{n=1}^{\infty}n
  t_n\xi^{-n-1} + O(1)\right)d\xi,
\]
and the $\Omega_k(P)$ are defined by
\begin{align}\begin{split}
    \Omega_k(P)&= \frac{\mathrm{d}\xi}{\xi^{k+1}}-
    \sum_{m=1}^{\infty}\frac{1}{m}\omega_{km}\xi^{m-1}\mathrm{d}\xi,\\
    &\oint_{\mathfrak{a}_j}\Omega_k(P)=0, \; j=1,\ldots,g. \end{split}
\end{align}

\begin{theorem} [\bf Bilinear Identity]
Let $S^1$ be a small circle about the point $p\in X$, then
\begin{eqnarray}
  &&\int_{S^1}
  \theta\left(\mathfrak{A}(Q) - \mathfrak{A}(P) - \sum_{i=1}^{\infty}
    \boldsymbol{U}_i
    t_i-\boldsymbol{u}\right)\theta\left(\mathfrak{A}(Q) -\mathfrak{A}(P)
    -\sum_{i=1}^{\infty}\boldsymbol{U}_i t_i'+\boldsymbol{u}\right) \nonumber\\
  &&\nonumber\\
  &&\qquad\qquad\qquad\qquad\times\mathrm{\exp}
  \left\{\int_{P_0}^P\Omega(\boldsymbol{t}+\boldsymbol{t}',P)\right\}
E^{-2}(P,Q) \mathrm{d}\xi(Q) =0,\label{bilinear}
\end{eqnarray}
where $\mathfrak{A}(P)=\int_{\infty}^P\mathrm{d}\boldsymbol{v}$ -
Abelian image of a point $P$, $\xi(Q)$ is the local coordinate of the
point $Q$ near $P$, $\xi(P)=0$, $E(P,Q)$ is the Schottky-Klein prime
form (\ref{schottky}) and the integration is around the unit circle
centred at the point $p$, $S^1=\{ q , |\zeta(q)|=1 \}$.

\end{theorem}

\noindent{\bf Proof}. Note that the integrand in (\ref{bilinear}) is
holomorphic in $q$ in the interior of $S^1$ except at its centre $Q=P$
where it has an isolated essential singularity with zero residue.
Therefore this integral vanishes.

To prove the equivalence of (\ref{bilinearresidue}) and
(\ref{bilinear}) note first that
\begin{equation}
\Omega(\boldsymbol{t}+\boldsymbol{t}')
=\Omega(\boldsymbol{t})+\Omega(\boldsymbol{t}').
\end{equation}
In (\ref{bilinear}), the factor
\begin{equation}
  \theta\left(\mathcal{A}(Q)-\mathcal{A}(P) -\sum_{i=1}^{\infty}
    \boldsymbol{U}_i t_i -\boldsymbol{u}\right)\mathrm{\exp}
  \left\{\int_{P_0}^P\Omega(\boldsymbol{t})\right\}
\end{equation}
can be expressed as
\begin{equation}\mathrm{exp}\left\{ \sum_{n=1}^{\infty} t_n \xi^{-n} \right\}
\theta\left(\mathcal{A}(Q)-\mathcal{A}(P)-\sum_{i=1}^{\infty}\boldsymbol{U}_i
  t_i -\boldsymbol{u}\right)\mathrm{\exp}
\left\{\sum_{m,n=1}^{\infty}
\omega_{mn} t_n \frac{\xi^m}{n} \right\}.\end{equation}
Taking into account the form (\ref{tau}) for $\tau(\boldsymbol{t};\boldsymbol{u})$,
the last two factors in the above formula may be expressed
\begin{equation}
\theta\left(\mathcal{A}(Q)-\mathcal{A}(P)-\sum_{i=1}^{\infty}\boldsymbol{U}_i
  t_i -\boldsymbol{u}\right)\mathrm{\exp}
\left\{\sum_{m,n=1}^{\infty}\omega_{mn} t_n \frac{\xi^m}{n} \right\}
=\mathrm{exp}\left\{-\sum_{n=1}^{\infty} \frac{\xi^n}{n}
\frac{\partial^n}{\partial t_n}
\right\}\tau(\boldsymbol{t};\boldsymbol{u}),
\end{equation}
which completes the proof of the equivalence of (\ref{bilinear}) and
(\ref{bilinearresidue}).

\begin{remark} Equations of this form were first written down
in terms of `vertex operators'
$$
\mathrm{exp}\left\{\sum_{1}^{\infty} t_n x^{-n}\right\}
\mathrm{exp}\left\{-\sum_{1}^{\infty}\frac{x^n}{n}
  \frac{\partial}{\partial t_n}\right\}
$$
for general $\tau$-functions
in the work of Sato. They may be understood as generating functions
for integrable hierarchies of Hirota bilinear equations.
\end{remark}

Using this theorem we may obtain partial differential equations
relating the Kleinian symbols, $\wp_{ij}$, $\wp_{ijk}$ etc. To do this
we now substitute the $\tau$-function (\ref{tausigma}) into this
expression and compute the residue.  This is parameterized by
$\boldsymbol{t}$ and $\boldsymbol{e}$.  The coefficients of monomials
in $\boldsymbol{t}$ are differential-difference expressions in
$\sigma$ at arguments $\boldsymbol{+u}$ and $\boldsymbol{-u}$.  To get
purely differential expressions, we now let $\boldsymbol{u}
\rightarrow 0$; the 'time' derivatives act on $\sigma$ via its
$\boldsymbol{u}$-dependence, though, so $\sigma_i(\boldsymbol{u}) =
-\sigma_i(\boldsymbol{-u})\rightarrow \pm \sigma_i(0)$, etc.  We then
replace the derivatives of the $\sigma$-function by derivatives of the
$\wp$ function using the recursive relations (\ref{sigmaladder})
described earlier.  This recovers the partial differential equations
for the Kleinian functions which we need.

{\sf Example I}

In the genus 2 hyperelliptic case, the first nonvanishing relations
occur at $t$-weight 3 in the $t_i$, i.e. the coefficients of $t_1^3$ or
$t_3$.  Both these give the first 4-index relation
\[
\wp_{1111} = 6\wp_{11}^2+\alpha_4 \wp_{11}+4\wp_{12}+\tfrac12\alpha_3.
\]
The next non-zero term is at weight 5, where we recover the relation
for $\wp_{1112}$, etc.  In contrast with the previous section, this
approach only gives the even derivative relations at odd $t$-weights.

The calculations by this approach, because they involve bilinear
products of $\tau$-functions, soon become very computer-intensive,
and we have not pursued them very far.

{\sf Example II}
Within this method we are able to recover all
quadratic relations for four indexed symbols $\wp_{ijkl}$, $1\geq
i,j,k,l\leq 3$ and cubic relations for three indexed symbols
$\wp_{ijk}^2$. We report here the {\em quartic} relation between even
variables, which must be one of the relations defining the Kummer
variety of the (3,4)-curve,
\begin{align*}
  \wp_{12}^{4}&- \wp_{22}^{3}-2\
  \wp_{11}\wp_{33}+2\wp_{13}^{2}+4\
  \wp_{12}^{2}\wp_{13}+\wp_{1113}\wp_{22}-6\wp_{11}
  \wp_{13}\wp_{22}+4\wp_{11}\wp_{12}\wp_{23}\\ & + \wp_{11}^{2}
  \wp_{22}^{2}-2\wp_{11}\wp_{12}^{2}\wp_{22}-\tfrac43\wp_{1113}
  \wp_{11}^{2}+8\wp_{11}^{3}\wp_{13}+2\mu_{12}+4\mu_{6}
  \wp_{11}^{3}-4\mu_{6}
  \wp_{11}\wp_{22}\\ & +3\mu_{3}\wp_{12}\wp_{13}+3\mu_{3}
  \wp_{11}\wp_{23}+\mu_{3} \wp_{12}^{3}+2\
  \mu_{6}\wp_{13}+\mu_{6} \wp_{12}^{2}+
  \mu_{9}\wp_{12}-{\mu_{3}}^{2}
  \wp_{11}^{3}\\ & -3\mu_{3}\wp_{11}\wp_{12}\wp_{22}=0.
\end{align*}
This relation is of independent interest since it is of weight 12 and
{\bf cannot} be written in the form (\ref{kum25}).  The detailed
structure of the Kummer variety in the (3,4) case is receiving further
investigation and will be reported on elsewhere.

\section{Discussion}
In view of the ``compare and contrast'' objectives of our paper, the
reader may wish to know which of these methods is the most effective
when calculating the required PDEs.  Unfortunately this is not an easy
question to answer.  Currently the two methods associated with the
$\tau$-function take rather longer to execute, with the residue method
in particular being slow and with large memory overheads.  However
this may be due in part to the fact that we have extensive experience
over ten years or more in developing computer algebra code for the
``Classical method'' in all its variations.  We have much less
experience in working with the tau-function methods.  So it may be
that with more study and with more computational experience, the
$\tau$-function methods become more competitive.  We would stress,
however, that both $\tau$-function methods are more systematic than
the classical methods, and may provide a suitable way of calculating
the PDEs associated with more complicated curves.

{\bf Acknowledgments} The work was partially supported by the European
Science Foundation Programme MISGAM (Methods of Integrable Systems,
Geometry and Applied Mathematics), the Imperial College (London)
section of the ENIGMA network, and a Fellowship (VZE) in the
Hanse-Wissenschaftskolleg in Delmenhorst in 2010, where the final
version of the paper was prepared.  The authors are grateful to L.\
Haine for discussions and for bringing to our attention the useful
unpublished manuscript by J.\ Fay, \cite{fay83}.  VZE would like to
thank John Harnad for many discussions and inspiring ideas, as well as
A.\ Nakayashiki for sending his paper on a closely related subject
before publication.  Some of this work was carried out whilst JCE and
VZE were visiting Iwate University, and these authors are grateful to
Y.\ Onishi for hospitality and for financial support from JSPS grant
19540002.


\bibliographystyle{plain}

\end{document}